\documentclass[10pt]{article}
\usepackage{amssymb,amsbsy,amsmath,amsfonts,amssymb,amsthm,mathrsfs}
\usepackage{latexsym,euscript,exscale}

\usepackage{times}

\newcommand {\ca} {{2^\N}}

\newcommand {\N}{\mathbb N}
\newcommand {\Q}{\mathbb Q}
\newcommand {\R}{\mathbb R}
\newcommand {\Z}{\mathbb Z}

\newcommand {\E}{\mathbb E}
\newcommand {\D}{\mathbb D}

\newcommand{\eps}{\epsilon}

\newcommand{\compl}{\complement}

\newcommand{\iso}{\cong}

\newcommand{\tom} {\emptyset}

\newcommand{\begr}{\!\upharpoonright}

\newcommand{\saa}{\Longrightarrow}
\newcommand{\equi}{\Longleftrightarrow}

\newcommand{\til}{\rightarrow}

\newcommand {\del}{ \; \big| \;}

\newcommand {\for}{\bigcup}

\newcommand {\snit}{\bigcap}

\newcommand {\ku} {\mathcal}

\newcommand{\ov}{\overline}
\newcommand{\inv}{^{-1}}

\renewcommand {\a} {\forall}

\newcommand{\conti}{2^{\aleph_0}}

\newcommand{\pf}{

\smallskip

\noindent {\it Proof : }}

\newcommand{\pff}{$\hfill  \Box$

\smallskip }

\newtheorem{thm}{Theorem}
\newtheorem{cor}[thm]{Corollary}
\newtheorem{lemme}[thm]{Lemma}
\newtheorem{prop} [thm] {Proposition}
\newtheorem{defi} [thm] {Definition}

\theoremstyle{remark}
\newtheorem{claim}{Claim}[thm]

\newtheorem*{cla}{Claim}

\begin{document}
\title{Automatic continuity of homomorphisms and fixed points on metric compacta}

\author{Christian Rosendal and S\l awomir Solecki
\footnote{Research of the second author was supported by NSF grant
DMS-0400931.}\\
University of Illinois at Urbana--Champaign\\
rosendal@math.uiuc.edu and ssolecki@math.uiuc.edu}
\date {January 2005}
\linespread {1.0}

\maketitle
\begin{abstract}
We prove that arbitrary homomorphisms from one of the groups ${\rm Homeo}(\ca)$, ${\rm Homeo}(\ca)^\N$, ${\rm Aut}(\Q,<)$, ${\rm Homeo}(\R)$, or ${\rm Homeo}(S^1)$ into a separable group are automatically continuous. This has consequences for the representations of these groups as discrete groups. For example, it follows, in combination with a result on V.G. Pestov, that any action of the discrete group ${\rm Homeo}_+(\R)$ by homeomorphisms on a compact metric space has a fixed point.
\end{abstract}

\section{Introduction}

The classical theorem of Pettis \cite[Theorem 9.10]{kec} says that
any \emph{Baire measurable} homomorphism from a Polish group to a
separable group is continuous. Measure theoretic counterparts of
this result are also known. Recently, it was proved \cite{sol}
that if $G$ is an amenable at $1$ Polish group, then any
\emph{universally measurable} homomorphism from $G$ to a separable
group is continuous. All locally compact Polish groups and Abelian
Polish groups are amenable at $1$ as are, for example, countable
products of amenable locally compact Polish groups. Therefore, the
classical measure theoretic automatic continuity results of Weil
(locally compact groups) and Christensen (Abelian groups) are
contained in this theorem. Going beyond homomorphisms with
regularity assumptions, of Baire category type or of measurable
type, in automatic continuity results requires the domain group to
be complicated. It was discovered recently \cite{kecros} that
\emph{each} homomorphism whose domain is a Polish group with ample
generic and whose range is separable is continuous. A Polish group
$G$ \emph{has ample generics} if for each finite $n\geq 1$ the
diagonal conjugation action of $G$ on $G^n$ given by
$$
g\cdot(h_1,\ldots,h_n)=(gh_1g\inv,\ldots,gh_ng\inv)
$$
has a comeagre orbit. Examples of groups with ample generics can
be found among permutation groups; most importantly $\ku
S_\infty$, the group of all permutations of $\N$, is such a group
(see \cite{kecros}).

The principal goal of this paper is to exhibit groups which are
less complicated than those with ample generics but which still
have this very strong automatic continuity property. However, our
results have broader consequences: for extremely amenable groups
and for questions concerning representing groups as subgroups of
$\ku S_\infty$ or as linear groups.

The main tool in our proofs of automatic continuity is a version
of the classical fact due to Steinhaus that if $A\subseteq \R$ is
a measurable set of positive Lebesgue measure, then $0\in{\rm
Int}(A-A)$. (Analogous lemmas were proved later by Pettis for
non-meager subsets $A$ of Polish groups with $A$ having the Baire
property and by Weil for non-Haar zero Haar measurable subsets of
locally compact groups.)

\begin{defi}Let $G$ be a topological group. We say that $G$ is
{\em Steinhaus} if there is $k\geq 1$ such that for any symmetric
countably syndetic set $W\subseteq G$ (i.e., covering $G$ by
countably many left-translates), $W^k$ contains an open
neighbourhood of $1_G$. To emphasise the exponent $k$, we will
sometimes say that $G$ is Steinhaus with exponent $k$.
\end{defi}

The first examples of Steinhaus groups come from \cite{kecros} where
it was shown that Polish groups with ample generics are Steinhaus
with exponent $10$.

The proposition below makes a connection between the Steinhaus
property and automatic continuity. Its proof is analogous to the
derivation of continuity of Baire measurable homomorphisms from
Pettis' lemma.

\begin{prop}\label{continuity} Let $G$ be a Steinhaus
topological group and $\pi:G\til H$ a homomorphism into a
separable group. Then $\pi$ is continuous.\end{prop}

\pf We need only show that $\pi$ is continuous at $1_G$. So
suppose $U\subseteq H$ is an open neighbourhood of $1_H$ and find
some symmetric open $V$ such that $1_H\in V\subseteq
V^{2k}\subseteq U\subseteq H$. As $H$ is separable, $V$ covers $H$
by countably many translates $\{h_nV\}$. So for each $h_nV$
intersecting $\pi[G]$ take some $g_n\in G$ such that $\pi(g_n)\in
h_nV$. Then $h_nV\subseteq \pi(g_n)V\inv V=\pi(g_n)V^2$, and hence
the $\pi(g_n)V^2$ cover $\pi[G]$. Now, if $g\in G$, find $n$ such
that $\pi(g)\in \pi(g_n)V^2$, whence $\pi(g_n\inv g)\in V^2$ and
$g_n\inv g\in\pi\inv[V^2]$, i.e., the $g_n\pi\inv[V^2]$ cover $G$.

So $W=\pi\inv[V^2]$ is symmetric and countably syndetic and hence
$1_G\in {\rm Int}(W^k)$ for some $k\geq 1$. But then
$\pi[W^{k}]\subseteq V^{2k}\subseteq U$, and therefore, $1_G\in
{\rm Int}(\pi\inv[U])$ and $\pi$ is continuous at $1_G$.\pff

We see that the exact condition we need to impose on $H$ in the
above proof is that any non-empty open set covers $H$ by countably
many translates. This condition is known in the literature as
$\aleph_0$-boundedness and is equivalent to embedding as a
subgroup into a direct product of second countable groups (see Guran \cite{gur}). Thus, by the definition of the product
topology, to show that any homomorphism with range in an
$\aleph_0$-bounded group is continuous, it is enough to show that
any homomorphism into a second countable group is continuous. So
from our perspective the three notions of second countable,
separable, and $\aleph_0$-bounded are equivalent.

Let us mention an immediate corollary.

\begin{cor}\label{C:pres} Let $G$ and $H$ be Polish groups. If $G$
is Steinhaus and $H$ is an image of $G$ by a homomorphism, then
$H$ is Steinhaus.
\end{cor}

\pf The homomorphism between $G$ and $H$ is continuous, by
Proposition~\ref{continuity}, and therefore, since it is surjective
and since $G$ and $H$ are Polish, it is open. Now it follows that
$H$ is Steinhaus by a straightforward computation. \pff

The main result of this paper is the following.

\begin{thm}\label{steinhausexamples} The following groups are
Steinhaus:
$$
{\rm Homeo}(\ca),\, {\rm Homeo}(\ca )^{\N},\, {\rm Aut}(\Q,<),\,
{\rm Homeo}_+(\R),\,\hbox{ and }{\rm Homeo}_+(S^1).
$$
\end{thm}

Here, ${\rm Homeo}_+(\R)$ is the group of increasing homeomorphisms
of $\R$, and similarly ${\rm Homeo}_+(S^1)$ is the group of
orientation preserving homeomorphisms of the unit circle, both of
them with the topology of uniform (or equivalently, pointwise)
convergence.

A crucial role in our arguments will be played by the existence of
comeager conjugacy classes in ${\rm Aut}(\Q ,<)$ \cite{tru2},
\cite{kustru}, and in ${\rm Homeo}(\ca )$ and ${\rm Homeo}_{+}(\R )$
\cite{kecros}. Of course, having ample generics implies having a
comeager conjugacy class, however, it is known that ${\rm Aut}(\Q
,<)$ and ${\rm Homeo}(\R )$ do not have ample generics; whether
${\rm Homeo}(\ca)$ has ample generics is open. All conjugacy classes
of ${\rm Homeo}_{+}(S^1)$ are meager.

\begin{cor}An arbitrary homomorphism from either ${\rm Homeo}(\ca)$, ${\rm Homeo}(\ca)^\N$,
${\rm Aut}(\Q,<)$, ${\rm Homeo}(\R)$, or ${\rm Homeo}(S^1)$ into a
separable group is continuous.\end{cor}

\pf The result is clear for the three first groups. And for the
last, notice that any homomorphism $\pi:{\rm Homeo}(\R)\til H$ into
a separable group restricts to a continuous homomorphism from
${\rm Homeo}_+(\R)$ into $H$, whence $\pi$ is also continuous as
${\rm Homeo}_+(\R)$ is open in ${\rm Homeo}(\R)$. Similarly for ${\rm
Hom}(S^1)$.\pff

Methods used to prove Theorem~\ref{steinhausexamples} yield also
the following result on the small index property.

\begin{thm}\label{smallindex}${\rm Homeo}_+(\R)$ is the only proper
subgroup of ${\rm Homeo}(\R)$ of index $<2^{\aleph_0}$. Similarly,
${\rm Homeo}_+(S^1)$ is the only proper subgroup of ${\rm Homeo}(S^1)$
of index $<2^{\aleph_0}$.\end{thm}

\section{Applications}

\subsection{Connections with extreme amenability}

A phenomenon that has recently received considerable attention in
topological dynamics is extreme amenability \cite{Pesbook}. A
topological group is called \emph{extremely amenable} if all of
its continuous actions on compact Hausdorff spaces have fixed
points. Such groups are also said to have the {\em fixed point on
compacta property}. The first examples of these groups date back
to work of Herer and Christensen \cite{hc}. Of special interest to
us are the results of Pestov \cite{pest} stating that ${\rm
Aut}(\Q,<)$ and ${\rm Homeo}_+(\R)$ are extremely amenable.

All these examples belong necessarily to the topological setting
as Ellis \cite{ellis} proved that any abstract (that is, discrete)
group acts freely on a compact space. (This was generalised later
by Veech \cite{vee} to locally compact groups.) Therefore,
non-trivial discrete groups are not extremely amenable. However,
using our results on automatic continuity, we will see that when
one restricts the attention to actions on metrisable compacta,
extreme amenability type phenomena occur for abstract groups.

\begin{cor} An arbitrary action by homeomorphisms of
${\rm Aut}(\Q,<)$ or ${\rm Homeo}_+(\R)$ on a compact metrisable
space has a fixed point.
\end{cor}

\pf It is enough to notice that any action by homeomorphisms on a
compact metric space $X$ corresponds to a homomorphism into the
homeomorphism group of $X$, which is Polish. By automatic
continuity the homomorphism is continuous, whence the action is
continuous. So, by the results of Pestov, there is a fixed point
on $X$.\pff

A continuous action of a Hausdorff topological group $G$ on
compact Hausdorff spaces is called a \emph{$G$-flow}. Such a flow
is called \emph{minimal} if each of its orbits is dense and it is a general fact that each flow
contains a minimal subflow. For any (Hausdorff) group $G$ there is a
{\em universal minimal flow}, that is, a continuous minimal action
of $G$ on a compact (Hausdorff) space $X$ such that any minimal
$G$-flow is a homomorphic continuous image of this one. So
extremely amenable groups are precisely those groups whose minimal
flows are one point spaces. There are Polish groups that are not
extremely amenable, but for which one nevertheless can compute the
universal minimal flow, and for some groups these flows turn out to be
metrisable. For example, Pestov showed that the universal minimal
flow of ${\rm Homeo}_+(S^1)$ is simply the canonical action of ${\rm
Homeo}_+(S^1)$ on $S^1$. Glasner and Weiss \cite{glawei} computed
the universal minimal flow of ${\rm Homeo}(\ca)$ to be equal to the
canonical action on the space of maximal chains of compact subsets
of $\ca$, a space which was first introduced and studied by
Uspenskij \cite{usp}. In general, there is no reason for a
topological Hausdorff group to have a universal minimal
\emph{metrisable} flow, that is, a metrisable flow which maps
homomorphically onto any metrisable flow. However, using the above results
of Pestov and Glasner--Weiss and our automatic continuity result,
we obtain the existence of such universal minimal metrisable flows
for groups whose universal minimal flows are non-metrisable.

\begin{cor} The discrete groups ${\rm Homeo}_+(S^1)$ and
${\rm Homeo}(\ca)$ have universal minimal \emph{metrisable} flows,
namely the canonical action on $S^1$ and the canonical action on
the space of maximal chains of compact subsets of $\ca$,
respectively.\end{cor}

Note that the universal minimal flows of the discrete groups
considered in the corollary above are non-metrisable, as  Kechris,
Pestov and Todorcevic proved in \cite{kpt} that the universal
minimal flow of any locally compact, non-compact group is
non-metrisable.

\subsection{Connections with representations}

There is some interest in representing groups by (faithful)
actions on countable sets which corresponds to embedding the
groups into $\ku S_\infty$. The maximal size of such groups is
$\conti$, and it is well-known that not every group of size at
most $\conti$ is embeddable into $\ku S_\infty$. For example, the
quotient of $\ku S_\infty$ by its subgroup of finitary
permutations is not so embeddable. Vladimir Pestov recently
informed us that John D. Dixon had asked for a characterisation of
the separable (Hausdorff) topological groups that were abstractly
embeddable into $\ku S_\infty$, and, in particular, he asked if
there was any counter-example at all. Theorem \ref{smallindex}
shows that this is so. For clearly, if ${\rm Homeo}_+(\R)$ acts on a
set $X$ of size $\kappa<\conti$, then the isotropy subgroup of
${\rm Homeo}_+(\R)$ at some $x\in X$ is of index $\leq \kappa$ and
hence has to be the whole group. Thus ${\rm Homeo}_+(\R)$ has no
non-trivial representations by permutations on a set of size
$<\conti$ whatsoever.

Consider now a more general type of representations namely linear
representations. Megrelishvili \cite{megrelishvili}, proving a
conjecture of Pestov, has shown that  ${\rm Homeo}_+(\R)$ has no
non-trivial strongly continuous representations by linear
isometries on a reflexive Banach space. Thus by Theorem
\ref{steinhausexamples}, we see that that ${\rm Homeo}_+(\R)$ has no
representations by linear isometries on a reflexive separable
Banach space.
\begin{cor}\label{isometr}
${\rm Homeo}_+(\R)$, as an abstract group, has  no non-trivial
representations by permutations on a set of size $<\conti$.
Moreover, it has no representations by linear isometries on a
reflexive separable Banach space.
\end{cor}

\subsection{Homomorphisms into locally compact Polish groups}

Pestov's theorem that there are no fixed point free continuous
actions of ${\rm Aut (\Q,<)}$ or ${\rm Homeo}_+(\R)$ on compact
spaces and Veech's theorem that each locally compact group acts
freely and continuously on a compact space in connection with our
automatic continuity result (and the remarks following
Proposition~\ref{continuity}) imply that there is no abstract
non-trivial homomorphism from these groups to locally compact
$\sigma$-compact groups. (Note that for locally compact groups
$\aleph_0$-boundedness coincides with $\sigma$-compactness.)
Similarly, a non-trivial homomorphism from a group $G$ into a
locally compact group $H$, induces a representation of $G$ by
linear isometries on $L^2$ of the Haar measure of $H$, which is
separable if $H$ is second countable. So the slightly weaker result
that ${\rm Homeo}_+(\R)$ has no nontrivial homomorphism into a
second countable locally compact group also follows from
Corollary~\ref{isometr}.

We show now that the same result holds for ${\rm Homeo}(\ca)$. In
fact, the argument for it is direct and applies also to ${\rm Aut
(\Q,<)}$ and ${\rm Homeo}_+(\R)$. We will use the following
property of an abstract group $F$ first studied by Bergman
\cite{ber}:

\begin{quote}Whenever $W_0\subseteq W_1\subseteq \ldots \subseteq F$
is an exhaustive sequence of subsets, then for some $n$ and $k$,
$F=W_n^k$.\end{quote}

\noindent That this condition holds for ${\rm Homeo}(\ca)$ is due
to Droste and G\"obel
\cite{drogob}. It also holds for ${\rm
Aut}(\Q,<)$ and ${\rm Homeo}(\R)$ as proved by Droste and Holland \cite{drohol}.

\begin{thm}\label{loc-comp}
Let $G$ be a Polish group which has the above property and has a
comeager conjugacy class. Then there is no non-trivial abstract
homomorphism from $G$ into a locally compact $\sigma$-compact
group.
\end{thm}

\pf Suppose $\pi: G\til H$ is a homomorphism into a locally
compact Polish group. We claim that $\ov{\pi[G]}$ is compact. To
see this take some increasing exhaustive sequence of compact
subsets of $H$:
$$
K_0\subseteq K_1\subseteq \ldots\subseteq H
$$
Then also
$$
(K_0\cap\pi[G]) \subseteq (K_1\cap\pi[G])\subseteq \ldots\subseteq
{\pi[G]}
$$
is exhaustive, so, since $G$ is Bergman, there are $n$ and $k$ such
that $\pi[G]=(K_n\cap\pi[G])^k\subseteq K_n^k$. But then $\pi[G]$
is relatively compact, whence $\ov{\pi[G]}$ is compact.

Secondly, we claim that the group $K=\ov{\pi[G]}$ has a dense
conjugacy class. For let $C\subseteq G$ be the comeagre conjugacy
class of $G$ and suppose $V\subseteq K$ is some non-empty open
set. We claim that $V\cap\pi[C]\neq \tom$. This suffices as
$\pi[C]$ is contained in a single conjugacy class of $K$. First
notice that as $K$ is compact and $\pi[G]$ is dense, there are
$\{g_i\}_{i\leq n} \subseteq G$ such that $K=\for_{i\leq
n}\pi(g_i)V$, whence $G=\for_{i\leq n}g_i\pi\inv(V)$. As $C$ is
comeagre, $\snit_{i\leq n}g_iC\neq\tom$, so take some
$h\in\snit_{i\leq n}g_iC\neq\tom$ and let $m$ be such that $h\in
g_m\pi\inv(V)$. Then clearly, $g_m\inv h\in \pi\inv(V)\cap C$, and
thus $\pi(g_m\inv h)\in V\cap \pi[C]$.

Therefore $\ov{\pi[G]}$ is compact with a dense conjugacy class.
But any conjucacy class in a compact Hausdorff group is closed, so
$\ov{\pi[G]}=\{1\}$ as the conjugacy class of $1$ is $\{1\}$.
\pff

The result above implies that groups with Bergman's property and
with a comeager conjugacy class cannot have non-trivial
representations by automorphisms of locally finite graphs.
Similarly, they cannot act non-trivially by isometries on compact
metric spaces.

\section{${\rm Homeo}(\ca)$ and ${\rm Homeo}(\ca )^{\N}$}

Before we begin our proofs, let us first mention the following
elementary fact, which will be used repeatedly.

\begin{lemme}\label{L:prod} Suppose  $n=1,2,\ldots, \aleph_0$ and
$\{G_i\}_{i\leq n}$ are Polish groups with comeagre conjugacy
classes. Then $G=\prod_{i\leq n} G_i$ has a comeagre conjugacy
class.\end{lemme}

\pf Let for each $i\leq n$, $\ku O_i\subseteq G_i$ be the comeagre
conjugacy class. Then obviously, $\prod_{i\leq n} \ku O_i\subseteq
G$ is a conjugacy class of $G$. Moreover, as
$$
\prod_{i\leq n} \ku O_i=\snit_{i\leq n}\Big[\big(\prod_{j\neq
i}G_j\big)\times \ku O_i\Big]
$$
by Kuratowski-Ulam this class is comeagre in $G$.\pff

\begin{thm}\label{hom}${\rm Homeo}(\ca)$ is Steinhaus.
\end{thm}

\pf Let $G={\rm Homeo}(\ca)$ and assume $W\subseteq G$ is symmetric
and covers $G$ by countably many left-translates $k_nW$, $k_n\in
G$. In particular, $W$ cannot be meagre and must therefore be
dense in some non-empty open set in $G$. But then $W\inv W=W^2$ is
dense in an open neighbourhood of the identity in $G$ and we can
therefore find some finite subalgebra $\mathsf A\subseteq {\rm
clopen}(\ca)$ with atoms $A_1,\ldots,A_p$ such that $W^2$ is dense
in $G_{(\mathsf A)}=\{g\in G\del g[A_i]=A_i, \;i=1,\ldots, p\}$.

For each $i=1,\ldots, p$ choose a point $x_i\in A_i$ and let
$B_n^i\subseteq A_i$ be a sequence of disjoint clopen sets
converging in the Hausdorff metric to the set $\{x_i\}$. Moreover,
let $B_n=B_n^1\cup\ldots\cup B_n^p$.

\begin{claim}\label{CL:ful} For some $n$, $B_n$ is full for
$W^2$, i.e., if $\gamma\in \{ g\in {\rm Homeo}(B_n): g[B^i_n] =
B^i_n, i=1, \dots p\}$, then there is a $g\in W^2$ such that
$g\begr B_n=\gamma$.
\end{claim}

\noindent Proof of Claim~\ref{CL:ful}. It will suffice to show
that some $B_n$ is full for $k_nW$ since then it is clearly also
full for $W^2 = (k_nW)^{-1}k_nW$. Assuming otherwise, we can for
each $n$, find some $\gamma_n\in{\rm Homeo}(B_n)$ such that for all
$g\in k_nW$, $g\begr B_n\neq\gamma_n$. Due to the convergence of
the sets $B_n$ to $\{x_1, \ldots, x_p\}$, there is a $f\in G$ such
that $f\begr B_n=\gamma_n$ and
$$
f\begr (\ca\setminus \for_n B_n)={\rm id}_{\ca\setminus \for_n
B_n}
$$
But then by the choice of $\gamma_n$, $f\notin k_nW$ for any $n$,
contradicting that the sets $k_nW$ cover $G$, and the claim is
proved.\pff

So we can choose some $B=B_{n_0}$ that is full for $W^2$. Let now
$G(B)=\{g\in G_{(\mathsf A)}\del {\rm supp}(g)\subseteq B\}$.

\begin{claim}\label{CL:more} $G(B)\subseteq W^{12}$.
\end{claim}

\noindent Proof of Claim~\ref{CL:more}. To see this, we notice
first that $G(B)$ is topologically isomorphic to ${\rm Homeo}(\ca)^p$ and, therefore by \cite{kecros} and
Lemma~\ref{L:prod}, it has a comeagre conjugacy class. Find now
some $n_1$ such that $k_{n_1}W\cap G(B)$ is non-meagre in $G(B)$.
Then

\begin{displaymath}\begin{split} (k_{n_1}W\cap G(B))\inv
\cdot(k_{n_1}W\cap G(B))&=(W\inv k_{n_1}\inv\cap
G(B))\cdot(k_{n_1}W\cap G(B))\\&\subseteq W^2\cap G(B)\end{split}
\end{displaymath}
is also non-meagre in $G(B)$, so we can find $f_0\in W^2\cap
G(B)$, whose conjugacy class in $G(B)$ is comeagre in $G(B)$.

Take any $h\in G(B)$ and, by fullness of $B$ for $W^2$ proved in
Claim~\ref{CL:ful}, find $g\in W^2$ such that $h\begr B=g\begr B$.
Then as $f_0\begr{\compl B}={\rm id}_{\compl B}$, we have
$$
hf_0h\inv =gf_0g\inv \in W^6
$$
This means that the conjugacy class of $f_0$ in $G(B)$ is
contained in $W^6$ and that therefore $W^6\cap G(B)$ is comeagre
in $G(B)$, whence, as the square of a comeagre set is everything, $W^{12}\cap G(B)=G(B)$, proving the claim.\pff

Since $W^2\cap G_{(A)}$ is dense in $G_{(A)}$, we can pick an
$h_0\in W^2\cap G_{(\mathsf A)}$ such that $h_0[B]=\compl B$,
i.e., such that $h_0[B\cap A_i]=\compl B\cap A_i$ for each
$i=1,\ldots , p$. Then clearly,
\begin{displaymath}\begin{split}
G(\compl B)&=\{ g\in G_{(\mathsf A)}\del {\rm supp}(g)\subseteq
\compl
B\}\\
&= h_0G(B)h_0\inv\\
&\subseteq W^{16}
\end{split}\end{displaymath}
and
$$
G_{(\mathsf B)}=G(B)\cdot G(\compl B)\subseteq W^{28}
$$
where $\mathsf B$ is the subalgebra of ${\rm clopen}(\ca)$
generated by $\mathsf A$ and the set $B$. So $W^{28}$ contains the
open neighbourhood of the identity, $G_{(\mathsf B)}$, and ${\rm
Homeo}(\ca)$ is Steinhaus with exponent $28$.\pff

We now show how techniques similar to the ones employed above give
that the group ${\rm Homeo}(\ca)^\N$ is Steinhaus. This result
implies the previous one by Corollary~\ref{C:pres}. Moreover,
notice that it is not in general true that the countable product
of Steinhaus groups is Steinhaus. The simplest counter-example is
$({\mathbb Z}/2)^\N$. Here an ultrafilter on $\N$ corresponds to a
subgroup of index $2$, which is open if and only if the
ultrafilter is principal. So $({\mathbb Z}/2)^\N$ is not
Steinhaus, even though the discrete group ${\mathbb Z}/2$ is.
However, it seems plausible that if $W$ is a symmetric countably
syndetic subset of a product $\Pi_iG_i$ of Steinhaus topological
groups of some common exponent $k$, then $W^k$ should contain a
product $\Pi_iU_i$ of open subsets of the $G_i$. We have neither a
proof nor a counter-example.

\begin{thm}\label{Hom} ${\rm Homeo}(\ca)^\N$ is Steinhaus.
\end{thm}

\pf The group ${\rm Homeo}(\ca )^\N$ is isomorphic to the subgroup
$G$ of ${\rm Homeo}(\ca \times\N)$ consisting of all $h\in {\rm
Homeo}(\ca\times \N)$ with $h[\ca\times \{ i\}] = \ca\times \{ i\}$
for each $i\in \N$. Put $K_i = \ca\times \{ i\}$. Let $W\subseteq
G$ be such that $G= \bigcup_nk_nW$ for some $k_n\in G$, $n\in \N$.

We will borrow two things from the proof of Theorem~\ref{hom}.
First note that the proof of Theorem~\ref{hom} gives that for any
$m\in \N$ a relatively open neighborhood of the identity of the
subgroup
$$
\{ f\in G\del \hbox{supp}(f)\subseteq \bigcup_{i< m}K_i\}
$$
is contained in $W^{28}$. For this reason, it will suffice to
prove that there exists $m$, perhaps depending on $W$ and the
sequence $(k_n)$, such that
$$
\{ f\in G\del \hbox{supp}(F)\subseteq \bigcup_{m\leq i}K_i\}\subseteq
W^{108}.
$$

Second, note that the following claim can be proved just like
Claim~\ref{CL:more} in the proof of Theorem~\ref{hom}. The only
difference is that in an appropriate place we need to use
Lemma~\ref{L:prod} with $n= \aleph_0$.

\begin{claim}\label{CL:fulseq} Let $U^n_i\subseteq K_n$ for
$i\leq n$ be pairwise disjoint clopen sets. There exists $n_0$
such that
$$
\{ f\in G\del {\rm supp}(f)\subseteq \bigcup_{n\geq n_0} U^n_{n_0}\}
\subseteq W^{12}.
$$
\end{claim}

Now, let ${\mathcal U}, {\mathcal U}'$ be families of pairwise disjoint
clopen subsets of $\bigcup_nK_n$. We say that ${\mathcal U}'$
\emph{refines} $\mathcal U$ if $\bigcup {\mathcal U} = \bigcup
{\mathcal U}'$ and each set in ${\mathcal U}'$ is included in a
set from $\mathcal U$. If ${\mathcal U}'$ refines $\mathcal U$ and
$\sigma$ and $\tau$ are permutations of $\mathcal U$ and
${\mathcal U}'$, respectively, we say that $\tau$ \emph{refines}
$\sigma$ if $\tau(U')\subseteq \sigma(U)$ whenever $U'\subseteq U$
for $U'\in {\mathcal U}'$ and $U\in {\mathcal U}$. Finally, if
$\mathcal U$ refines $\{ K_i\del i\geq n\}$ for some $n\in \N$, let
${\rm Sym}_0({\mathcal U})$ be the group of all permutations of
$\mathcal U$ which refine the identity permutation of $\{ K_i\del i\geq n\}$.

\begin{claim}\label{CL:diag} There exist $n_0\in\N$ and a family of
clopen sets $\mathcal U$ refining $\{ K_i\del i\geq n_0\}$ such that
for any ${\mathcal U}'$ refining $\mathcal U$ and any $\tau \in
{\rm Sym}_0({\mathcal U}')$ refining ${\rm id}\in {\rm
Sym}_0({\mathcal U})$ we can find $h\in W^2$ with
$$
h[U] = \tau(U)\;\hbox{ for any }U\in {\mathcal U}'.
$$
\end{claim}

\noindent Proof of Claim~\ref{CL:diag}. It suffices to find
$n_0\in\N$, a family of clopen sets $\mathcal U$ refining $\{ K_i\del
i\geq n_0\}$, and $\sigma\in {\rm Sym}_0({\mathcal U})$ such that
for any ${\mathcal U}'$ refining $\mathcal U$ and any $\tau \in
{\rm Sym}_0({\mathcal U}')$ refining $\sigma$ we can find $h\in
k_{n_0}W$ with
$$
h[U] = \tau(U)\;\hbox{ for any }U\in {\mathcal U}'.
$$
Indeed, the claim follows from the statement above since
$$
W^2 = (k_{n_0}W)^{-1} k_{n_0}W\;\hbox{ and }\;\sigma^{-1}\circ
\sigma = {\rm id}\in {\rm Sym}_0({\mathcal U}).
$$

Assume towards a contradiction that the statement fails. Let
${\mathcal U}_{0}$ be $\{ K_i\del  i\geq 0\}$ and let $\tau_{0}\in
{\rm Sym}_0({\mathcal U}_{0})$ be the identity. Assume we are
given ${\mathcal U}_{n}$ refining $\{ K_i\del i\geq n\}$ and
$\tau_n\in {\rm Sym}({\mathcal U}_n)$ such that there is no $h\in
k_{n-1}W$ with $h[U] = \tau(U)$ for all $U\in {\mathcal U}_n$.
Consider
$$
{\mathcal U}'_{n} = \{ U\in {\mathcal U}_n\del U\subseteq
\bigcup_{i\geq n+1} K_i\}\; \hbox{ and }\; \tau_n' =
\tau_n\upharpoonright {\mathcal U}'_{n}.
$$
By our assumption, we can find ${\mathcal U}_{n+1}$ refining
${\mathcal U}'_{n}$ and $\tau_{n+1}\in {\rm Sym}_0({\mathcal
U}_{n+1})$ refining $\tau_n'$ such that for no $h\in k_{n+1}W$ do
we have $h[U] = \tau(U)$ for all $U\in {\mathcal U}_{n+1}$.

The inductive construction allows us to find $h_0\in G$ such that
for each $U\in {\mathcal U}_n$ with $U\subseteq K_n$ we have
$h_0[U] = \tau_n(U)$. Note that, again by the construction,
$h_0\not\in k_nW$ for each $n$. This yields a contradiction since
$\bigcup_n k_nW = G$, and the claim follows.

\begin{claim}\label{CL:five} Let $B, C\subseteq 2^\N$ be clopen.
Assume that $C\cap B$ and $C\setminus B$ are both non-empty. Let
$G_1 =\{ f\in {\rm Homeo}(\ca)\del f[B]=B\}$ and $G_2 = \{ f\in {\rm
Homeo}(\ca )\del {\rm supp}(f)\subseteq C\}$. Then ${\rm Homeo}(\ca ) =
G_1G_2G_1G_2G_1$.
\end{claim}

\noindent Proof of Claim~\ref{CL:five}. Let $f\in {\rm Homeo}(\ca)$.
Consider the clopen sets $L_1 = f[B]\cap B$ and $L_2 =
f[\complement B]\cap \complement B$. Assume both of them are
non-empty. This assumption allows us to find $g_1\in G_1$ such
that
$$
g_1[B\setminus L_1] \varsubsetneq C\cap B\;\hbox{ and }\;
g_1[\complement B\setminus L_2]\varsubsetneq C\cap \complement B.
$$
Now there exists $g_2\in G_2$ such that
$$
g_2[g_1[B\setminus L_1]] \subseteq C\cap \complement B\;\hbox{ and
} \; g_2[g_1[\complement B\setminus L_2]]\subseteq C\cap B
$$
and
$$
g_2[C\setminus g_1[B\setminus L_1]] \subseteq C\cap B\;\hbox{ and
} \; g_2[C\setminus g_1[\complement B\setminus L_2]]\subseteq
C\cap \complement B.
$$
Note that $f(g_2g_1)^{-1}\in G_1$, so the conclusion follows.

Assume now $L_1$ or $L_2$ is empty, say $L_1=\emptyset$. Then note
that $\complement B \setminus L_2\not=\emptyset$. Let $g_1\in G_1$
be such that
$$
g_1[\complement B \setminus L_2]\cap C\not=\emptyset.
$$
If now $g_2\in G_2$ is such that $g_2[C\setminus B] = C\cap B$,
then clearly the sets $L_1$ and $L_2$ computed for
$f(g_2g_1)^{-1}$ are both non-empty, so we can apply the previous
procedure to get the conclusion of the claim.

We prove now the theorem from the three claims. Pick $n_0$ and
$\mathcal U$ as in Claim~\ref{CL:diag}. For $n\geq n_0$, pick
pairwise disjoint non-empty clopen sets $V^n_i\subseteq K_n$ with
$i\leq n$ in such a way that for each $U\in {\mathcal U}$ with
$U\subseteq K_n$ and each $i\leq n$ we have $U\cap V^n_i
\not=\emptyset$. We assume no $V^n_i$ contains a $U\in {\mathcal
U}$. By Claim~\ref{CL:fulseq}, there is $n_1\geq n_0$ such that
$$
\{ f\in G\del \hbox{supp}(f)\subseteq \bigcup_{n\geq n_1}
V^n_{n_1}\}\subseteq W^{12}.
$$
For $U\in {\mathcal U}$ with $U\subseteq K_n$ for some $n\geq
n_1$, put
$$
U^0 = V^n_{n_1}\cap U\;\hbox{ and }\; U^1 = U\setminus V^n_{n_1}.
$$
Let $B_0 =\bigcup \{ U^0\del U\in {\mathcal U}\hbox{ and } U\subseteq
\bigcup_{n\geq n_1} K_n\}$ and $B_1 =\bigcup \{ U^1\del U\in
{\mathcal U}\hbox{ and } U\subseteq \bigcup_{n\geq n_1} K_n\}$. So
we have
\begin{equation}\label{E:12}
\{ f\in G\del \hbox{supp}(f)\subseteq B_0\} \subseteq W^{12}.
\end{equation}
Note that ${\mathcal U}' = \{ U^0, U^1\del U\in {\mathcal U}\}$
refines $\mathcal U$ and $\tau \in {\rm Sym}_0({\mathcal U}')$
given by $\tau(U^j) = U^{1-j}$ refines $\rm id$. Therefore, by
Claim~\ref{CL:diag}, we have $h_0\in W^2$ such that $h_0[U^j] =
\tau(U^j)$ for each $U^j\in {\mathcal U}'$. It follows that
$h_0[B_0] = B_1$ and $h_0[B_1] = B_0$, whence from \eqref{E:12}
\begin{equation}\label{E:16}
\begin{split}
\{ f\in G\del \hbox{supp}(f)\subseteq B_1\}\,&\subseteq h_0^{-1}\{
f\in G\del \hbox{supp}(f)\subseteq B_0\}h_0\\
&\subseteq W^2W^{12}W^2 =W^{16}.
\end{split}
\end{equation}

For $n\geq n_1$ pick pairwise disjoint clopen sets $C^n_i\subseteq
K_n$ for $i\leq n$ so that each $C^n_i$ intersects both $B_0$ and
$B_1$. Applying Claim~\ref{CL:fulseq}, we see that there is
$n_2\geq n_1$ such that
\begin{equation}\label{E:77}
\{ f\in G\del  \hbox{supp}(f)\subseteq \bigcup_{n\geq n_2} C^n_{n_2}\}
\subseteq W^{12}.
\end{equation}
Now using Claim~\ref{CL:five} for each $n\geq n_2$ (with $B=
B_0\cap K_n$ and $C= C^n_{n_2}$) along with \eqref{E:12},
\eqref{E:16}, and \eqref{E:77}, we get
\begin{equation}\notag
\{ f\in G\del \hbox{supp}(f)\subseteq \bigcup_{n\geq n_2} K_n\}
\subseteq W^{28}W^{12}W^{28}W^{12}W^{28} = W^{108},
\end{equation}
and the theorem follows.\pff

Let us now see how this result leads to automatic continuity for
other groups containing ${\rm Homeo}(\ca)^\N$. Fix a
denumerable model-theoretical structure $\mathsf A$ and suppose
that ${\rm Aut}(\mathsf A)$ is Steinhaus, or just that any
homomorphism from ${\rm Aut}(\mathsf A)$ into a separable group is
continuous. We can assume that the domain of $\mathsf A$ is $\N$.
Now let $\alpha$ be an action of ${\rm Aut}(\mathsf A)$ on ${\rm
Homeo}(2^\N)^\N$ defined as follows: $\alpha(g, [n\mapsto
h_n])=[n\mapsto h_{g(n)}]$. Thus, we can form the topological
semidirect product ${\rm Aut}(\mathsf A)\ltimes_\alpha {\rm
Homeo}(\ca)^\N$. Recall that the topology on the semidirect product
is the same as the product topology on ${\rm Aut}(\mathsf A)\times
{\rm Homeo}(\ca)^\N$. Now, if $K$ is a separable group and $\pi:{\rm
Aut}(\mathsf A)\ltimes_\alpha{\rm Homeo}(\ca)^\N\til K$ is a
homomorphism, then $\pi$ restricts to a continuous homomorphism on
each of the factors, whence $\pi$ is continuous on the semi-direct
product. When $\mathsf A$ is just the empty structure we have:

\begin{cor}
Let $X$ be the topological space $\N\times \ca$ and $\mathsf E$
the equivalence relation on $X$ given by $(n,\alpha)\mathsf
E(m,\beta)\equi n=m$. Then any homomorphism from ${\rm Homeo}(X,\mathsf E)$
(i.e., the group of homeomorphisms of $X$ preserving the
equivalence relation $\mathsf E$) into a separable group is
continuous.
\end{cor}

\section{${\rm Aut}(\Q,<)$}

\begin{thm}\label{Q}${\rm Aut}(\Q,<)$ is Steinhaus.
\end{thm}

Our proof of this result relies on the combinatorics of Truss' proof
from \cite{tru1} that ${\rm Aut}(\Q,<)$ satisfies the so called
small index property, that is, that every subgroup of index
strictly less than the continuum is open. Let $G={\rm Aut}(\Q,<)$
and let $\D$ be the family of all subsets $X\subseteq \Q$ of the
form
$$
X=\for_{n\in\Z}]x_{2n},x_{2n+1}[
$$
where $(x_n)_{n\in \Z}$ is a sequence of irrationals satisfying
$x_n<x_{n+1}$ and $x_n\til \pm\infty$ for $n\til\pm\infty$.
Moreover, for $X\in \D$ we let
$$
A(X)=\{g\in G\del {\rm supp}(g)\subseteq X\}
$$
Since any element $g\in {\rm Aut}(\Q,<)$ can be extended to a
unique homeomorphism of $\R$, we will sometimes evaluate
expressions $g(x)$ for $g\in G$ and $x$ an irrational number.

The following lemma can be extracted from \cite{tru1}.

\begin{lemme}\label{truss}(Truss)
$$
G=\for_{X,Y\in \D}A(X)\cdot A(Y)
$$
\end{lemme}

\pf Given $g\in G$, find a sequence $(x_n)_{n\in \Z}$ of
irrationals such that $x_{n-1}<g(x_n)<x_{n+1}$ and $x_n\til
\pm\infty$ for $n\til\pm\infty$. Now, put $I_n=]x_n,x_{n+1}[$ and
notice by the choice of $x_n$ that
$$
g\big([x_{4n},x_{4n+1}]\big)\subseteq ]x_{4n-1},x_{4n+2}[
$$
So we can define some $h\in G$ such that for each $n\in \Z$
$$
h\begr\; ]x_{4n+2},x_{4n+3}[={\rm id}
$$
$$
h\begr g(I_{4n})=g\inv
$$
Then $hg\begr I_{4n}={\rm id}$ and $h\inv\begr I_{4n+2}={\rm id}$.
Letting
$$
Y=\for_{n\in \Z}I_{4n+1}\cup I_{4n+2}\cup I_{4n+3}
$$
and
$$
X=\for_{n\in \Z}I_{4n}\cup I_{4n+1}\cup I_{4n+3}
$$
we have $g=h\inv\cdot hg\in A(X)\cdot A(Y)$.\pff

\

\noindent{\em Proof of Theorem \ref{Q}:} Suppose $W\subseteq G$ is
symmetric and countably syndetic. Then $W$ cannot be meagre and
hence $W^2$ must be dense in some open neighbourhood of the
identity, $V=\{g\in G\del g(q_1)=q_1,\ldots,g(q_{p-1})=q_{p-1}\}$,
for some rational numbers $q_1<\ldots<q_{p-1}$. Notice that $V$ is
topologically isomorphic to $G^{p}$. Fix some $k_n\in G$ such that
the $G= \bigcup_nk_nW$.

We now let $\E$ be the family of all sets $X\subseteq \Q$ of the
form
$$
X=\Big(\for_{n\in
\Z}]x^1_{2n},x^1_{2n+1}[\Big)\cup\ldots\cup\Big(\for_{n\in
\Z}]x^p_{2n},x^p_{2n+1}[\Big),
$$
where $x_n^i$ are irrationals such that for $q_0=-\infty$ and
$q_p=+\infty$, we have
$$
q_{i-1}<x^i_n<x^i_{n+1}<q_i,
$$
$$
x^i_n\til q_i \:\:{\rm for}\:\; n\til+\infty,
$$
$$
x_n^i\til q_{i-1} \;\;{\rm for} \;\;n\til -\infty.
$$
Moreover, for $X\in \E$, we let
$$
A(X)=\{g\in G\del {\rm supp}(g)\subseteq X\}.
$$
Clearly, by Lemma~\ref{truss}, we have
$$
V=\for_{X,Y\in \E}A(X)\cdot A(Y).
$$
So to prove that ${\rm Aut}({\Q})$ is Steinhaus with exponent
$96$, it suffices to show the following claim.

\begin{cla} $A(X)\subseteq W^{48}$ for any $X\in \E$.
\end{cla}

\noindent{Proof of Claim.} Fix some $X\in \E$ and sequences
$x_n^i$ as above. Moreover, for each $i=1,\ldots,p$, let
$$
I^i_n=]x^i_{2n},x_{2n+1}^i[
$$
and for each $\vec a=(a_1,\ldots, a_p)$, where $a_i\subseteq \Z$
is bi-infinite, let
$$
X_{\vec a}=\big(\for_{n\in
a_1}I^1_n\big)\cup\ldots\cup\big(\for_{n\in a_p}I^p_n\big).
$$

We stress the fact that the $I^i_n$ name only every second
interval of the $\Z$-ordered partition of $]q_{i-1},q_i[$ into the
intervals $]x_m^i,x^i_{m+1}[$. Thus, if $h\in A(X_{\vec a})$, then
$h[I^i_n] = I^i_n$ for each $n\in a_i$, $i= 1, \dots , p$.  Now
pick a sequence of $\vec a_n$ such that the sets $X_{\vec a_n}$
are all disjoint, which is equivalent to demanding that the $j$-th
terms of $\vec a_n$ and $\vec a_m$ are disjoint for $n\neq m$ and
$j=1,\ldots,p$. From the remark about $h$ above it follows that if
$g_n\in A(X_{\vec a_n})$, then $g: {\Q}\to {\Q}$ defined to be
$g_n\begr X_{\vec a_n}$ on $X_{\vec a_n}$ and the identity on
${\Q}\setminus \bigcup_n X_{\vec a_n}$ is an element of $G= {\rm
Aut}({\Q})$.

We claim that for some $n_0$, $X_{\vec a_{n_0}}$ is full for
$k_{n_0}W$, i.e., that for any $g\in A(X_{\vec a_{n_0}})$, there
is some $h\in k_{n_0}W$ such that $g\begr X_{\vec a_{n_0}}= h\begr
X_{\vec a_{n_0}}$. If not, we could for each $n$ find some $g_n\in
A(X_{\vec a_n})$ such that for all $h\in k_nW$, we have $g_n\begr
X_{\vec a_n}\neq h\begr X_{\vec a_n}$. As noticed above, we could
then find one single $g\in G$ such that $g\begr X_{\vec
a_n}=g_n\begr X_{\vec a_n}$ for every $n$. But this would
contradict that the $k_nW$ cover $G$.

So suppose $X_{\vec a_{n_0}}$ is full for $k_{n_0}W$. Then $X_{\vec
a_{n_0}}$ is also full for $W^2$. For simplicity, let $\vec
a=(a_1,\ldots, a_p)=\vec a_{n_0}$.

Clearly, $A(X_{\vec a})$ is topologically isomorphic to
$(G^\Z)^p$, so it has a comeagre conjugacy class by \cite{tru2}
and Lemma~\ref{L:prod}. Find now some $n_1$ such that
$k_{n_1}W\cap A(X_{\vec a})$ is non-meagre in $A(X_{\vec a})$,
whence also $W^2\cap A(X_{\vec a})$ is non-meagre in $A(X_{\vec
a})$. Therefore there is some $f\in W^2$ belonging to the comeagre
conjugacy class in $A(X_{\vec a})$. But if $h\in A(X_{\vec a})$,
then there is a $g\in W^2$ agreeing with $h$ on $X_{\vec a}$,
whence $hfh\inv=gfg\inv\in W^6$. So $W^6$ contains the comeagre
conjugacy class of $A(X_{\vec a})$ and as the product of two
comeagre sets in a group is everything, $A(X_{\vec a})\subseteq
W^{12}$.

Let now $(\vec a^\alpha)_\alpha$ be a continuum size family of
sequences $\vec a^\alpha=(a_1^\alpha,\ldots,a_p^\alpha)$ such that
$a_i^\alpha\subseteq a_i$ is bi-infinite and $a_i^\alpha\cap
a_i^\beta$ is finite for every $\alpha\neq \beta$. (see, e.g., Kunen \cite{kunen}, p. 48.)

For each $\alpha$ write also $\Q\setminus \{q_1,\ldots,q_{p-1}\}$
as a disjoint union of non-empty irrational intervals
$J^\alpha_{i,n}$ ($n\in \Z$, $i=0,\ldots, p$), such that
$$
J^\alpha_{0,n}<J^\alpha_{0,n+1}<q_1<J^\alpha_{1,n}<J^\alpha_{1,n+1}<
q_2<\ldots<q_{p-1}<J^\alpha_{p,n}<J^\alpha_{p,n+1},
$$
where
$$
X_{\vec a^\alpha}=\for_{\substack{i=0,\ldots,p\\
n\in\Z}}J^\alpha_{i,2n},
$$
$$
\Q\setminus X_{\vec a^\alpha}=\for_{\substack{i=0,\ldots,p\\
n\in\Z}}J^\alpha_{i,2n+1}.
$$
We notice that this forces each $J^\alpha_{i,2n}$ to be equal to
some $I^i_m=]x^i_{2m},x^i_{2m+1}[$ for an $m\in a^\alpha_i$, while
each $J^\alpha_{i,2n+1}$ must be on the form
$J^\alpha_{i,2n+1}=]x^i_{2m+1},x^i_{2l}[$ for some $m< l$ in
$a^\alpha_i$.

Now, find $g_\alpha\in V$ such that
$$
g_\alpha\big[J^\alpha_{i,n}\big]=J^\alpha_{i,n+1}
$$
By the uncountability there is $n_2$ and distinct $\alpha$ and
$\beta$ such that $g_\alpha, g_\beta\in k_{n_2}W$, whence
$g_\alpha\inv g_\beta, g_\beta\inv g_\alpha\in W^2$.

If $n\notin a^\alpha_i$, then $I^i_n\subseteq J^\alpha_{i,2l+1}$
for some $l$, whence
$$
g_\alpha\big[I^i_n\big]\subseteq
g_\alpha\big[J^\alpha_{i,2l+1}\big]=J^\alpha_{i,2l+2}=I^i_m
$$
for some $m>n$ with $m\in a^\alpha_i$. Similarly, if $n\notin
a^\beta_i$, then $g_\beta\big[I^i_n\big]\subseteq I^i_m$ for some
$m>n$ with $m\in a^\beta_i$. This, along with the almost
disjointness of $a^\alpha_i$ and $a^\beta_i$, allows us to find $N$ big enough so that for all
$i=0,\ldots, p$ and $|n|\geq N$
$$
n\notin a^\alpha_i \saa g_\alpha\big[I^i_n\big]\subseteq
I^i_m\;\;\; (m\in a^\alpha_i\setminus a_i^\beta)
$$
$$
n\notin a^\beta_i \saa g_\beta\big[I^i_n\big]\subseteq I^i_m\;\;\;
(m\in a^\beta_i\setminus a_i^\alpha)
$$
From this, by a similar argument, we get
$$
n\notin a^\alpha_i \saa g_\beta\inv
g_\alpha\big[I^i_n\big]\subseteq I^i_l\;\;\;(l\in
a_i^\beta\subseteq a_i)
$$
$$
n\notin a^\beta_i \saa g_\alpha\inv
g_\beta\big[I^i_n\big]\subseteq I^i_l\;\;\;(l\in
a^\alpha_i\subseteq a_i)
$$
Suppose also that $N$ has been chosen large enough to ensure that
for all $|n|\geq N$ either $n\notin a_i^\alpha$ or $n\notin
a_i^\beta$. Then for all $|n|\geq N$, either
\begin{equation}\label{1}
g_\beta\inv g_\alpha\big[I^i_n\big]\subseteq X_{\vec a}
\end{equation}
or
\begin{equation}\label{2}
g_\alpha\inv g_\beta\big[I^i_n\big]\subseteq X_{\vec a}
\end{equation}
As $W^2$ is dense in $V$, we can choose $h\in W^2$ and $m_i\in
a_i$ such that
$$
h\big[I^i_{-N}\cup\ldots\cup I^i_N\big]\subseteq I^i_{m_i}
$$
for every $i=0,\ldots,p$. Then for all $n\in\Z$, either \eqref{1}
or \eqref{2} or
\begin{equation}\label{3}
h\big[I^i_n\big]\subseteq X_{\vec a}
\end{equation}
where as noticed $g_\beta\inv g_\alpha,\; g_\alpha\inv g_\beta,\;
h\in W^2$ and $A(X_{\vec a})\subseteq W^{12}$.

Now define sets
\begin{align*}
b_i&=\{n\in \Z\del g_\beta\inv g_\alpha\big[I^i_n\big]\subseteq
X_{\vec a}\}\\
c_i&=\{n\in \Z\del g_\alpha\inv g_\beta\big[I^i_n\big]\subseteq
X_{\vec a}\}\\
d_i&=\{n\in \Z\del h\big[I^i_n\big]\subseteq X_{\vec a}\}
\end{align*}
and let $\vec b=(b_0,\ldots,b_p)$, $\vec c=(c_0,\ldots, c_p)$ and
$\vec d=(d_0,\ldots, d_p)$. Since \eqref{1}, \eqref{2}, or
\eqref{3} holds for each integer $n$, we get
$$
A(X)=A(X_{\vec b})\cdot A(X_{\vec c})\cdot A(X_{\vec d}).
$$
Since additionally, directly from the definitions of $b_i$, $c_i$
and $d_i$, we have
\begin{align*}
A(X_{\vec b})&\subseteq (g_\beta\inv g_\alpha)\inv A(X_{\vec
a})g_\beta\inv g_\alpha\\
A(X_{\vec c})&\subseteq (g_\alpha\inv g_\beta)\inv A(X_{\vec
a})g_\alpha\inv g_\beta\\
A(X_{\vec d})&\subseteq h\inv A(X_{\vec a})h,
\end{align*}
we get $A(X)\subseteq (W^{16})^3=W^{48}$, proving the claim and
thereby the theorem.\pff

\section{${\rm Homeo}(\R)$ and ${\rm Homeo}(S^1)$}

\begin{thm}\label{homeo}${\rm Homeo}_+(\R)$ is Steinhaus.
\end{thm}

\pf Let us first recall the Polish group topology on ${\rm
Homeo}(\R)$. It has as basis the following sets
$$
U(f;q_1,\ldots,q_{p-1};\eps)=\{g\in {\rm Homeo}(\R)\del
d(f(q_i),g(q_i))<\eps, \a i<p\}
$$
where $f\in {\rm Homeo}(\R)$, $\eps>0$, and $q_1<\ldots<q_{p-1}$
belong to $\R$. A similar topology on ${\rm Homeo}([0,1])$ gives a
topologically isomorphic group. The structure of the subgroup
${\rm Homeo}_+(\R)$ of all increasing homeomorphisms is very similar
to that of ${\rm Aut}(\Q,<)$ except from the fact that the former
is connected and the latter is totally disconnected. Nevertheless,
the proof for ${\rm Aut}(\Q,<)$ translates almost word for word
into a proof for ${\rm Homeo}_+(\R)$. Let us just mention the
changes needed:

The first thing to notice is that there is a comeagre conjugacy
class in ${\rm Homeo}_+(\R)$, which is shown in \cite{kecros}.
Secondly, instead of working with irrational intervals of $\Q$ one
replaces these by say half open intervals $]r,s]\subseteq \R$. One
easily sees that Lemma \ref{truss} goes through as before.
Supposing now that $W\subseteq {\rm Homeo}_+(\R)$ is symmetric and
countably syndetic, we find some open neighbourhood of the identity
$$
V=\{g\in {\rm Homeo}_+(\R)\del d(q_i,g(q_i))<\eps, \a i< p\}
$$
in which $W^2$ is dense. We can suppose that
$$
\eps< \min_i\frac{d(q_i,q_{i+1})}{3}
$$
Let also
$$
U=\{g\in {\rm Homeo}_+(\R)\del g(q_i)=q_i, \a i< p\}
$$
and notice that $U$ is topologically isomorphic to ${\rm
Homeo}_+(\R)^{p}$. One sees that
$$
U=\for_{X,Y\in \E}A(X)\cdot A(Y)
$$
where $\E$ is defined as in the proof of Theorem \ref{Q}. We can
now prove that $A(X)\subseteq W^{48}$ exactly as Claim was
established in the proof of Theorem~\ref{Q}, noticing that we do
not need $W^2$ to be dense in $U$ but only in $V$. So we get that
$U\subseteq W^{96}$. But $U$ is not open in the Polish topology on
${\rm Homeo}_+(\R)$. The following claim will show that ${\rm
Homeo}_+(\R)$ is Steinhaus with exponent $194$.

\begin{claim}\label{dingo}$V\subseteq W^{194}$\end{claim}

\noindent Proof of Claim~\ref{dingo} Suppose $f\in V$ and
$f(q_i)=r_i$ for $i=1,\ldots, p-1$. Then by the density of $W^2$
in $V$ there is an $h\in W^2$ such that
$$
d(r_i,h(q_i))<1/2d(r_i,q_i)
$$
for $i=1,\ldots, p$. But then there is also some $g\in U\subseteq
W^{96}$ satisfying $g(h(q_i))=r_i$, whence $f=(gh)h\inv g\inv
f\subseteq UW^2U\subseteq W^{194}$, since $h\inv g\inv f\in U$.
\pff

Consider now the group of homeomorphisms of the unit circle ${\rm
Homeo}(S^1)$ with the topology of uniform (or equivalently,
pointwise) convergence. As in the case of ${\rm Homeo}(\R)$, we
let ${\rm Homeo}_+(S^1)$ be the index $2$ subgroup of orientation
preserving homeomorphisms. It is a well-known fact that ${\rm
Homeo}_+(S^1)$ does not even have a non-meagre conjugacy class, as
e.g. the rotation number (see Katok and Hasselblatt, \cite[Chapter
11]{kathas}) is a continuous conjugacy invariant. But, on the
other hand, for any $x\in S^1$ the closed subgroup
$$
{\rm Homeo}(S^1,x)=\{g\in {\rm Homeo}(S^1)\del g(x)=x\}
$$
is topologically isomorphic to ${\rm Homeo}(\R)$.

\begin{thm}${\rm Homeo}_+(S^1)$ is Steinhaus.
\end{thm}

\pf We will deduce this result from the result for ${\rm
Homeo}_+(\R)$. So assume that $W\subseteq{\rm Homeo}_+(S^1)$ is
symmetric and countably syndetic. Then we can find some
neighbourhood of the identity
$$
V=\{g\in {\rm Homeo}_+(S^1)\del d(q_i,g(q_i))<\eps, \a i\leq n\}
$$
in which $W^2$ is dense. Let
$$
H=\{g\in {\rm Homeo}_+(S^1)\del g(q_1)=q_1\}
$$
then $H$ is topologically isomorphic to ${\rm Homeo}_+(\R)$.

\begin{claim}\label{unnamed} $W^2\cap H$ is symmetric and
countably syndetic in $H$.\end{claim}

\noindent Proof of Claim~\ref{unnamed}. Let $\{k_nW\}_\N$ cover
${\rm Homeo}_+(S^1)$. For each $n$ such that $k_nW\cap H\neq \tom$
pick some $g_n$ in the intersection. Then $k_n\in g_nW\inv =g_nW$
and thus
\begin{equation}\begin{split}
H&\subseteq \big(\for_n k_nW\big)\cap H\\
&\subseteq \big(\for_ng_nW^2\big)\cap H\\
&=\for_ng_n\Big(W^2\cap H\Big)
\end{split}\end{equation}
where the last equality holds as $g_n\in H$. The claim is proved.

Now, as $H\iso {\rm Homeo}_+(\R)$ is Steinhaus with exponent $194$,
there are $p_1,\ldots,p_m\in S^1$ and $\eps>\delta>0$ such that
\begin{equation}\notag
\begin{split}
U=\{g\in {\rm Homeo}_+(S^1)\del d(p_i,g(p_i))<\delta, \a i\leq m
\;\&\; g(q_1)=q_1\}\,&\subseteq (W^2)^{194}\\
&=W^{388}.
\end{split}
\end{equation}
In particular,
$$
U'=\{g\in {\rm Homeo}_+(S^1)\del \a\leq n\; g(q_i)=q_i\;\&\; \a
i\leq m\; g(p_i)=p_i\}\subseteq W^{388}
$$
and $W^2$ is dense in the set
$$
V'=\{g\in {\rm Homeo}_+(S^1)\del \a\leq n\;
d(g(q_i),q_i)<\delta\;\&\; \a i\leq m\; d(g(p_i),p_i)<\delta\}
$$
As in the proof of Claim~\ref{dingo}, we see that $V'\subseteq
U'W^2U'\subseteq W^{768}$, hence ${\rm Homeo}_{+}(S^1)$ is Steinhaus
with exponent $768$. \pff

\begin{thm}${\rm Homeo}_+(\R)$ is the only proper subgroup of ${\rm
Homeo}(\R)$ of index $<2^{\aleph_0}$.\end{thm}

\pf Since ${\rm Homeo}_+(\R)$ is connected, it is enough to show
that any subgroup $H\leq {\rm Homeo}_+(\R)$ of index $<2^{\aleph_0}$
is open. This is done by repeating the proof above for $W=H$.
However, there are a few things that have to be noticed before
this can be done. Namely, since $H$ is not necessarily
countably syndetic, we have to see exactly where this is used and
propose a substitute. First of all, we have to prove that $H$
cannot be meagre in ${\rm Homeo}_+(\R)$, and secondly, we have to
show that some $X_{\vec a}$ is full for $H^2=H$.

To see that $H$ is not meagre in ${\rm Homeo}_+(\R)$, notice that
the map $(g,h)\mapsto gh^{-1}$ is continuous and open from $({\rm
Homeo}_+(\R))^2$ to ${\rm Homeo}_+(\R)$. So if $H$ is meagre then
$E=\{(g,h)\in ({\rm Homeo}_+(\R))^2: gh^{-1}\in H\}$ would be a
meagre equivalence relation and therefore by Mycielski's Theorem
(see \cite[Theorem 19.1]{kec}) have a continuum of classes,
contradicting $[{\rm Homeo}_+(\R):H]<2^{\aleph_0}$.

Now to see that some $X_{\vec a}$ is full for $H^2=H$, we pick our
sequence $\vec a_n$ as in the proof of Claim of Theorem~\ref{Q}
such that the sets $X_{\vec a_n}$ are all disjoint. Let now
$$
N=\{g\in{\rm Homeo}_+(\R)\del \a n\;g[X_{\vec a_n}]=X_{\vec
a_n}\}\leq {\rm Homeo}_+(\R)
$$
and notice that $N$ is topologically isomorphic to $\prod_n
A(X_{\vec a_n})$, which is itself isomorphic to ${\rm
Homeo}_+(\R)^\N$. Let now $H_n$ be the projection of $H\cap N$ into
$A(X_{\vec a_n})$. Then $H\cap N\leq\prod_n H_n$, whence
$$
\prod_n[A(X_{\vec a_n}):H_n]=[N:\prod_n H_n]\leq [N:H\cap N]\leq
[{\rm Homeo}_+(\R):H]<2^{\aleph_0}
$$
Therefore at most finitely many $[A(X_{\vec a_n}):H_n]$ can be
different from $1$, meaning that at least for some $n$,
$H_n=A(X_{\vec a_n})$, i.e., $X_{\vec a_n}$ is full for $H$.\pff

\begin{cor}${\rm Homeo}_+(S^1)$ is the only proper subgroup of ${\rm
Hom}(S^1)$ of index $<2^{\aleph_0}$.\end{cor}

\pf Fix two points, e.g. $i$ and $-i$, on the unit circle $S^1$
and suppose $H<{\rm Homeo}(S^1)$, $[{\rm Homeo}(S^1),H]<\conti$. For
$x\in S^1$ let
$$
{\rm Homeo}_+(S^1,x)=\{g\in {\rm Homeo}_+(S^1)\del g(x)=x\}
$$
which is a subgroup isomorphic to ${\rm Homeo}_+(\R)$. Moreover, we
have
$$
[{\rm Homeo}_+(S^1,i):H\cap {\rm Homeo}_+(S^1,i)]<\conti
$$
and
$$
[{\rm Homeo}_+(S^1,-i):H\cap {\rm Homeo}_+(S^1,-i)]<\conti
$$
Thus as ${\rm Homeo}_+(\R)$ has no proper subgroups of index
$<\conti$, ${\rm Homeo}_+(S^1,i)\leq H$ and ${\rm Homeo}_+(S^1,-i)\leq
H$. But it is not hard to see that
$$
{\rm Homeo}_+(S^1)={\rm Homeo}_+(S^1,i)\cdot{\rm
Homeo}_+(S^1,-i)\cdot{\rm Homeo}_+(S^1,i).
$$
So ${\rm Homeo}_+(S^1)=H$. \pff

\end{document}